\input amstex
\input epsf
\documentstyle{amsppt}
\NoBlackBoxes

\def \De {\Delta}

\def \rP {{\Bbb R} P}
\def\ls{\leqslant}
\def\gs{\geqslant}
\def \cP {{\Bbb C} P}

\def\fy{\varphi}

\def\MRat{\operatorname{MRat}}
\def\Aut{\operatorname{Aut}}
\def\ee{\varepsilon}
\def \bR{\Bbb R}

\def\L{\Cal L}
\def\P{\Cal P}
\def\N{\Cal N}

\topmatter

\title {Counting real rational functions with all real critical values}
\endtitle
\author   Boris Shapiro*\; and Alek Vainshtein\dag
\endauthor
\affil
   * Department~of  Mathematics, University~of Stockholm,\\ 
S-10691, Sweden, {\tt shapiro\@matematik.su.se}\\
\dag  Departments of Mathematics and of Computer Science,
University of Haifa,  Israel 31905,
{\tt alek\@mathcs.haifa.ac.il}
\endaffil

\leftheadtext{Boris Shapiro and Alek Vainshtein}

\abstract
We study the number $\sharp_{n}^\bR$ of real rational degree $n$
functions (considered up to a linear fractional transformation of the
independent
variable) with a given set of $2n-2$ distinct real critical values.
We present a combinatorial interpretation of these numbers and pose
a number of related questions.
\endabstract

\endtopmatter

\document

\heading {\S 1. Introduction }\endheading

This paper continues the series \cite{A2--A6, S, SV} dealing with the
topology of the sets of ordinary, trigonometric, and exponential
polynomials with the maximal possible number of real critical values.
Recall that a degree $n$ rational function is any holomorphic degree
$n$ map $\cP^1\to\cP^1$. A rational function $f$ is called {\it real\/}
if both the image and the preimage are equipped with an
antiholomorphic involution, which is preserved  by $f$.
In analogy with \cite{A2--A4}, we call a real rational function
with the maximal possible number of real critical values a
   {\it rational M-function.}
We say that a rational M-function is  {\it generic\/}
if all its critical values are distinct. Thus, any degree $n$ generic
rational $M$-function has exactly $2n-2$ real and distinct critical values.
We call two rational functions $f_{1}$ and $f_{2}$ {\it  equivalent\/}
if $f_{1}={f_{2}}\circ L$, where $L$ is an arbitrary linear
fractional transformation of the independent variable. In the
present paper we study the number $\sharp_{n}^\bR$ of nonequivalent
generic degree $n$ rational $M$-functions with a given set of $2n-2$ real
distinct critical values. The following result is a prototype for the
main theorem of this note.

\proclaim{Theorem A {\rm (see~\cite{A3})}}
The number  $\sharp^{pol}_n$ of nonequivalent real degree $n$ polynomials
with the same set of $n-1$ real distinct critical values
   equals the number   $K_{n-1}$  of up-down
permutations of length $n-1$. These numbers are called
{\it Euler-Bernoulli numbers} and their generating function is
$$
\sum_{n=0}^\infty K_n\frac{t^n}{n!}=\sec t+ \tan t.
$$
\endproclaim

{\smc Remark.} Note the if we drop the condition that the polynomials are
real, then the number of nonequivalent  degree $n$ polynomials
with the same set of $n-1$  distinct critical values
equals $n^{n-3}$ and is one of the simplest Hurwitz numbers, see \cite{L}.
   Thus, even if all the critical values are real, only an
exponentially small fraction of
   all classes of complex polynomials with these critical values
   contains real representatives.

Denote the set of all generic rational
M-functions of degree $n$ by $\MRat_n$.
A  {\it  planar  chord diagram\/} of order $2m$
   is a circle $\Cal C$ embedded in $\Bbb R^2$ with fixed $2m$ points
   $\{v_{i}\}$ called
{\it vertices\/}  connected pairwise by $m$ nonintersecting curve
segments $\{\chi_{j}\}$   called {\it chords\/}.  The domain $\Cal D$  bounded
by $\Cal C$ is called {\it the basic disk\/}. The $2m$ arcs
$\{b_{k}\}$ in $\Cal C\setminus \bigcup_{j=1}^{2m}v_{j}$ are
called {\it boundary arcs\/}. The $m+1$ connected components
$\{f_{l}\}$ of $\Cal D\setminus
\bigcup_{i=1}^m \chi_{i}$  are called
{\it faces}. An {\it  automorphism\/} of a planar chord diagram  is a
homeomorphism of its basic disk sending vertices to vertices and chords
to chords.
An alternative way to represent planar chord diagrams is provided by
planar trees.
The  {\it associated planar tree\/}  of a given
planar chord diagram is the tree whose vertices are in
1-1-correspondence with the faces of the diagram, and are connected
if and only if the corresponding faces are adjacent.

A planar chord diagram is
called {\it properly oriented}, or {\it directed}, if all its
chords and boundary arcs are provided with an
orientation in such a way that the boundary of each face becomes
a directed cycle. Obviously, in order to direct a planar chord diagram
it suffices to direct any one of its edges. Therefore, there exist
exactly two possible ways of directing a diagram, which are opposite to
each other, i.e.~the second one is obtained from the first one by
reversing the direction of every edge. Once and for all fixing the
standard orientation of the plane, we call a face of a directed planar
chord diagram {\it positive\/} if the face lies to the left when we
traverse its boundary according to the chosen direction, and {\it
negative\/} otherwise.

   Consider a cyclically ordered labeling set $\Cal S=\{1\prec 2
\prec 3 \prec \ldots
   \prec 2m\prec 1\}$. A properly oriented
   planar chord diagram is called {\it properly labeled\/} if its $2m$
   vertices are labeled by pairwise distinct elements from $\Cal S$ so
that the labels
   arising on the boundary of each face traversed according to its
   orientation form a cyclically ordered subset of $\Cal S$, see Fig
1. Finally,
   two properly oriented and properly labeled planar chord diagrams are called
   {\it equivalent\/} if there exists a homeomorphism of the basic
   disks (in general, not preserving the orientation of the preimage
   disk) sending vertices to  vertices, chords to chords, and
preserving labels and orientations.

   {\smc Remark.} The group $C_{2m}$ of cyclic shifts of the labels in
   the set $\Cal S$ acts on the equivalence classes of properly oriented and
   labeled planar diagrams of order $2m$ (preserving the class of the
underlying properly oriented diagram). Note that this action is free if one
   forbids automorphisms of the underlying diagram.

   \vskip 15pt
\centerline{\hbox{\epsfysize=3,5cm\epsfbox{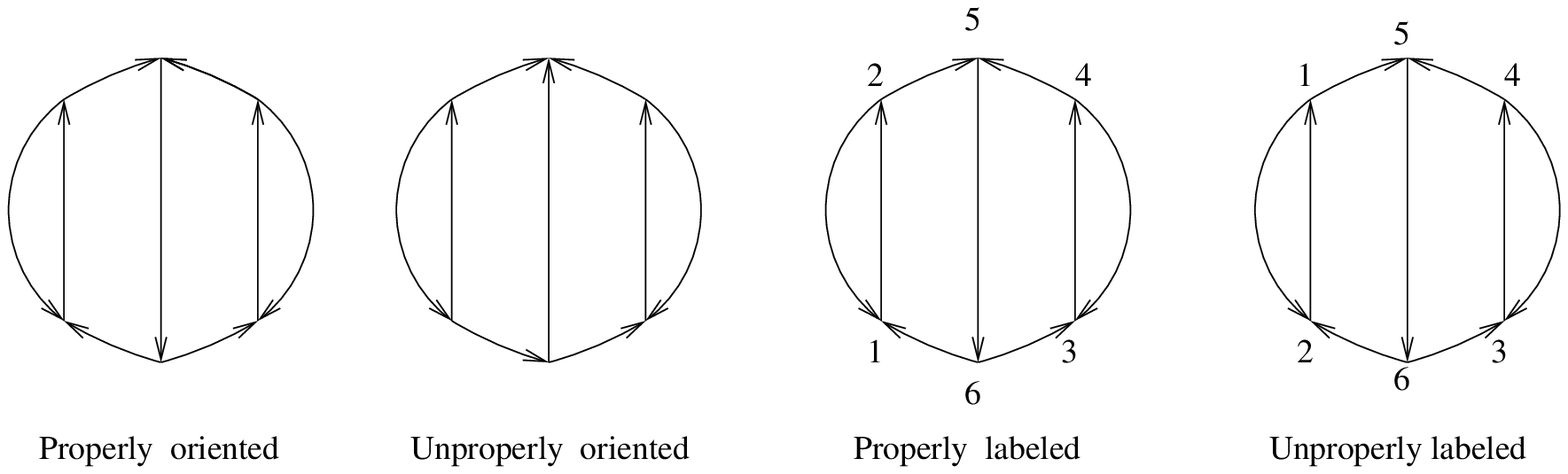}}}
\midspace{0.1mm} \caption{Fig.~1. Properly and unproperly oriented,
   properly and unproperly labeled planar diagrams }

Now we can formulate the main results of the present note.

\proclaim{Theorem 1} For $n \gs 3$
the number $\sharp_{n}^\bR$ of nonequivalent
generic degree $n$ rational M-functions with a given set of $2n-2$ real
distinct critical values equals the number of nonequivalent
properly oriented and labeled planar chord diagrams of order $2n-2$.
(The cases $n=1,2$ are exceptional and
$\sharp^{\bR}_{1}=\sharp^{\bR}_{2}=1$.)
\endproclaim

{\smc Remark.} Note that if we drop the condition that rational
functions are
real then the number of nonequivalent  degree $n$ rational functions
with the same set of $2n-2$  distinct critical values
   equals $n^{n-3}(2n-2)!/n!$ and is one of the simplest Hurwitz
   numbers, see e.g. \cite {CrT}.
   Thus even if all the critical values are real it seems plausible that
   only the exponentially small part of
   all classes of complex rational with these critical values contains
   real representatives.
\medskip

The numbers $\sharp^{\bR}_{n}$ can be considered as  real analogs
of the classical Hurwitz numbers, but they seem to be more difficult
to calculate. The first interesting values of $\sharp^{\bR}_{n},\,$
$n=3,4,5$ are $2,20,406$ resp.~(see examples below).

Theorem 1 presents the number  $\sharp_{n}^\bR$ as the sum of the
number of proper labelings over the set of all properly oriented
planar chord diagrams of order $2n-2$. The set of all properly
oriented  planar chord diagrams (considered up to plane
homeomorphisms) is rather an unconvenient object. There is a more
convenient base of summation for  $\sharp_{n}^\bR$ described below.

A planar chord diagram $D$ is said to be {\it rooted\/} if
one of its vertices is distinguished; the latter is then called the {\it root}.
A chord of $D$ is called a {\it diameter\/} if there exists a nontrivial
automorphism
of $D$ that fixes both endpoints of this chord; evidently, any diagram has at 
most one diameter.
Given a rooted planar chord diagram $D$, let $\sharp_r(D)$ denote the number 
of 
proper labelings of $D$ with label $1$ placed at the root and one of the two
possible orientations $ D^+$ or $D^-$ fixed.
The action of  the cyclic group on the set of all labels shows
that $\sharp_r(D)$ is
independent of the choice of orientation and of the root, provided the root
is not an endpoint of a diameter (in the latter case the number of proper 
labelings is  twice smaller). In other words, $\sharp_r(D)$ depends only on 
the planar chord diagram,
or equivalently, on the associated planar tree $T(D)$.
For this reason instead of $\sharp_r(D)$ we can write
$\sharp(T)$, where $T=T(D)$.

\proclaim{Theorem $\bold 1'$} For $n\gs3$ the
 number $\sharp_{n}^\bR$ of nonequivalent
generic degree $n$ rational $M$-functions with a given set of $2n-2$ real
distinct critical values equals
$$
\sharp_{n}^\bR=\sum_{T\in \P_{n}}\sharp(T),
$$
where $\P_n$ is the set of all planted trees on $n$ vertices.
\endproclaim

The number of connected components of the set $\MRat_n$ is given by the 
following result.

\proclaim{Theorem 2}
For $n\gs3$ the number $\widetilde\sharp^{\bR}_{n}$ of connected 
components of the set
$\MRat_{n}$ of all generic degree $n$ rational M-functions equals
the number of $C_{2n-2}$-orbits on the set of nonequivalent
properly oriented and labeled planar chord diagrams of order $2n-2$.
\endproclaim

The first interesting values of $\widetilde\sharp^{\bR}_{n},\,$
$n=3,4,5$ are $1,4,55$ resp.~(see examples below).

{\smc Remark.} The above Theorems~1 and~2 can be extended to a more
general situation of generic rational functions whose critical values
are not necessary real, see \cite {NSV}.
\medskip

{\smc Examples.} We start from the first nontrivial case $n=3$. In
this case there is only one planar chord diagram; its associated
planar tree contains two edges. This diagram is {\it
orientation-symmetric\/} (that is, it has an automorphism sending one
proper orientation to the other). It has 2 distinct proper labelings
that belong to the same unique $C_4$-orbit. So, there are two
different rational M-functions of degree 3, and the space $\MRat_3$ is
connected.

For $n=4$ there exist two
different planar chord diagrams whose trees are given in the right-hand
side of Figure~2. Both diagrams are orientation-symmetric.
The first of them has 18 proper labelings split into three full
$C_{6}$-orbits (i.e.
of length $6$) presented on the left upper part of Fig.2. The second one
has 2 proper labelings forming a single $C_{6}$-orbit. The total
number of distinct real M-functions of degree~4 is thus 20, and the
number of connected components in $\MRat_4$ equals~4.

   \vskip 15pt
\centerline{\hbox{\epsfysize=6.5cm\epsfbox{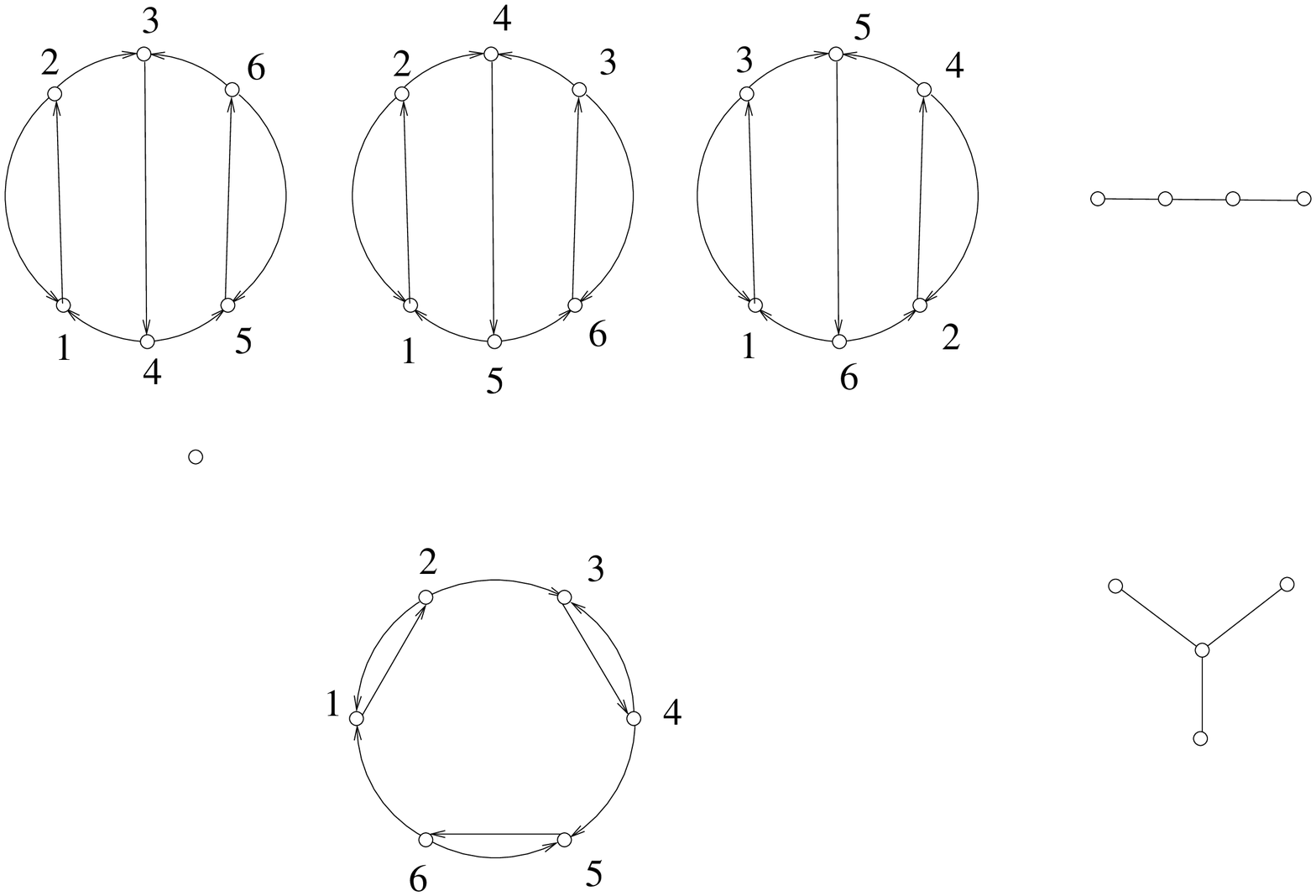}}}
\midspace{0.1mm} \caption{Fig.~2. Four $C_{6}$-orbits presenting
  $20$ properly labeled planar diagrams for $n=4$ }
\medskip

For $n=5$ there are 3 planar chord diagrams: the one with three square
faces, the one with a suqare face and a hexagonal face, and the one with an
octagonal face.
The first and the third of them are orientation-symmetric
while the second is not. The first diagram has $284$ proper labeling split
into $32$ full
orbits of $C_{8}$ (i.e. of length $8$) and seven $C_{8}$-orbits of
length $4$. The second diagram has
$120$ proper labelings split into $15$ full orbits of $C_{8}$.
Finally, the third diagram has two proper labelings forming a single
$C_{8}$-orbit. Therefore, there are 406 distinct real M-functions, and
the
number of connected components in $\MRat_5$ equals 55.

Another way to get the same number is to use Theorem~$1'$. There exist
3 planar trees on 5 vertices, which can be planted in 14 different ways,
see Fig.~3. The numbers $\sharp(T)$ for the three trees are 71, 15, and 1.
We thus get the same numbers 284, 120, and 2.

   \vskip 15pt
\centerline{\hbox{\epsfysize=5.5cm\epsfbox{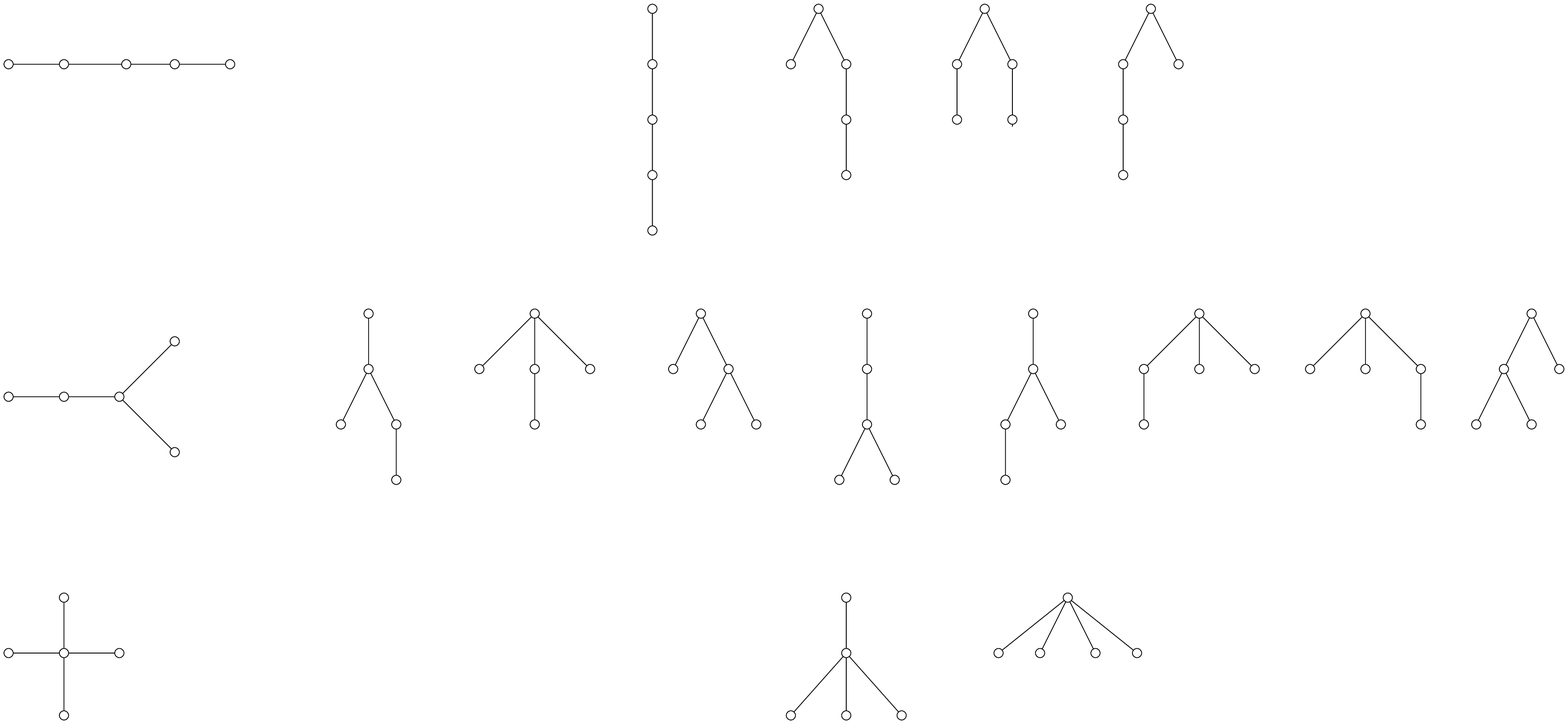}}}
\midspace{0.1mm} \caption{Fig.~3. Three planar trees on 5 vertices
 and their plantings }
\medskip

Below we present certain partial results on counting the number of
proper labelings of a given planar chord diagram.

We say that a planar chord diagram is of {\it $P_{k}$-type\/}
if its associated planar tree is a path on $k+1$ vertices. Let $\sharp(P_n)$
denote the number of nonequivalent properly oriented and labeled
planar chord diagrams of $P_n$-type.

Given a
permutation $\sigma\in S_n$ we define its {\it up-down sequence\/}
  as a word of length $n-1$ in the alphabet $\{U,D\}$ whose
$i$th letter  equals $U$ if and only if $\sigma_i<\sigma_{i+1}$.
We say that $\sigma\in S_{2n+1}$ is a {\it 2up-2down\/} permutation if
its up-down sequence equals $(UUDD)^k$ for $n=2k$ or $(UUDD)^kUU$ for
n=$2k+1$.

\proclaim {Theorem 3} For any $n\gs1$ the number $\sharp(P_n)$ equals
$n$ times the number of 2up-2down permutations of length $2n-1$.
\endproclaim

2up-2down and similar permutations were studied in \cite{CaS1, CaS2, CGJN}.
Theorem~3 together with results obtained in \cite{CaS2} yields
the following corollary. Denote by $F_{P}(x)$ the
modified exponential generating function
$$
F_{P}(x)=\sum_{n=1}^\infty\sharp(P_n)\frac{x^{2n}}{(2n)!}.
$$

\proclaim{Corollary 1} The generating function  $F_{P}(x)$ is given by
$$
F_{P}(x)=\frac x2\cdot\frac{\fy_0(x)\fy_1(x)-\fy_2(x)\fy_3(x)+\fy_3(x)}
{\fy_0^2(x)-\fy_1(x)\fy_3(x)},
$$
where
$$
\fy_j(x)=\sum_{k=0}^\infty \frac{x^{4k+j}}{(4k+j)!}\qquad j=1,2,3,4,
$$
is the $j$th Olivier function.
\endproclaim

Combining Corollary~1 with results obtained in \cite{CaS1} one gets the
following asymptotic estimate.

\proclaim{Corollary 2} $\sharp(P_n)\sim\alpha\frac{n(2n-1)!}{\gamma^{2n}}$,
where $\alpha$ is a constant and $\gamma=1.8750\dots$ is the smallest 
positive solution
of the equation $\cos z\cosh z+1=0$.
\endproclaim

More generally, consider planar chord diagrams whose
associated planar tree is
a caterpillar, that is, consists of a path  and of
an arbitrary number of edges incident to internal vertices of the path.
Let us order the vertices in the path of a caterpillar from one
endpoint to the other, and let $d_1,\dots,d_{k}$ be the sequence of
the degrees of the internal vertices of the caterpillar;
evidently, $d_i\gs2$ for $1\ls i\ls k$. 
In this case we say that the
  chord diagram is of type $(\delta_1,\dots,\delta_{k})$, where
$\delta_i= 2(d_i-1)$. 
In particular, for the $P_{k}$-type one has $\delta_i\equiv2$. 

Let $D$ be a planar chord diagram of type $ (\delta_1,\dots,\delta_{k})$,
and let $D^+$, $D^-$ be the two properly oriented diagrams corresponding to
$D$. Denote by $\sharp(D)$ the number of nonequivalent proper labelings of
$D^+$ and $D^-$.
The following result is a generalization of Theorem~3 for caterpillars.

\proclaim{Theorem 4} 
For any $k\gs1$,  the number
$\sharp(D)$ equals $n^*$ times the number of
permutations of length
$2e-1$ whose up-down sequence equals
$U^{\delta_1}D^{\delta_2}\dots$, where 
$$
n^*=\left\{\alignedat2 &\frac{2e}{|\Aut(D^+)|}\qquad & &
\text{if $D$ is orientation--symmetric,}\\
&\frac{4e}{|\Aut(D^+)|}\qquad & &
\text{otherwise,}
\endalignedat\right.
$$
$\Aut(D^+)$ is the group of automorphisms of the oriented diagram $D^+$,
and $e=1+\frac12\sum_{i=1}^{k}\delta_i$ is the number of edges in the 
corresponding caterpillar.
\endproclaim

Permutations with a given up-down sequence were studied in many
papers, starting with the classic work of MacMahon (see \cite{MM}).
His approach  leads to determinantal
formulas for the number of such permutations, rediscovered later by 
Niven \cite{Ni} from very basic
combinatorial considerations. For the relations of this approach to the
representation theory of the symmetric group see \cite{Fo, St}. Another, 
purely combinatorial 
approach to the same problem was suggested by
Carlitz \cite{Ca}. However, his general recursive formula for 
for the number of permutations with a given up-down sequence
is rather difficult to use.

Using the results of MacMahon and Niven,
we can reformulate Theorem~4 as follows.
Let 
$\ee_i=\delta_1+\delta_2+\dots+\delta_{i}$, and let 
$\{s_1<s_2<\dots\}$ denote the sequence
$$
\ee_1+1,\ee_1+2,\dots,\ee_2,\ee_3+1,
\dots, \ee_4,\dots, \ee_{2i-1}+1,\dots,\ee_{2i}, \ee_{2i+1}+1,
\dots;
$$
extend this
sequence by adding $s_0=0$ and by appending $2e-1$.

\proclaim{Corollary 3} For any $k\gs1$ and any diagram $D$ of type
$(\delta_1,\dots,\delta_k)$ one has
$$
\sharp(D)=n^*\cdot\det\left(\binom{s_i}{s_{j-1}}\right).
$$
\endproclaim

The authors are sincerely grateful to the Max-Planck Institut
f\"ur Mathematik in
Bonn for the financial support and the excellent research atmosphere
in the fall of 2000 when this project was initiated and in the summer
of 2002 when it was finished.
Sincere thanks goes to  A.~Eremenko who
suggested this  project to the first author in the fall of 1999.
The authors are also grateful to M.~Shapiro for numerous discussions.

\heading \S 2.  {General constructions and proofs} \endheading

\demo {Proof of Theorem 1}
   Let  $f\: \cP^1\to\cP^1$ be a real rational
function with all distinct real critical values (a generic {\it
M-function\/}).
Consider the set  $C(f)$ of its critical values.

Following \cite {Vi}, we define the  {\it net\/} $\N(f)$ of the function $f$
as the inverse image of the real cycle $\rP^1\subset \cP^1$.
Evidently,  $\N(f)$ contains $\rP^1$ and is invariant with respect to
  the complex
conjugation. Moreover, all the critical points of $f$ are real as
well; they belong to
$\rP^1\subset \N(f)$ and are the branching points of
   $\N(f)$ over $\rP^1$.
   Each critical point is adjacent to exactly $4$ arcs (edges)
belonging to $\N(f)$. These arcs  do not intersect outside
$\rP^1\subset \N(f)$
and split into pairs invariant under complex conjugation. Let us
fix some orientation of $\rP^1\subset \cP^1$ in the image. Then
all edges of $\N(f)$ will obtain the induced orientation.

Let us associate to a net  $\N(f)$ the properly oriented  and labeled
planar chord
diagram  $X(f)$ defined as follows. Roughly speaking, the chord
diagram is the part of $\N(f)$ contained in the closed upper halfplane.
More exactly, let us cyclically label the critical values along $\rP^1$
by the numbers $1,2,\ldots,2n-2$ from the labeling set $S$.
The number of vertices in  $X(f)$ equals $2n-2$ and they are in a
$1-1$-correspondence with the critical points of $\N(f)$. Two vertices
are connected by a chord if and only if the corresponding critical points
are connected by an arc belonging to $\N(f)$. The orientation of any
chord (including the boundary arcs) coincides with the induced
orientation of the arcs in $\N(f)$. Finally,
   each vertex is labeled by the label of the corresponding critical
   value. The connected components in  $\cP^1\setminus \N(f)$ are called the
{\it faces\/} of the net $\N(f)$.
Obviously, the faces of the net are split into complex conjugate
pairs, and each such pair corresponds to a single face in  $X(f)$.
Using a result by Caratheodory on the correspondence of the boundaries
under conformal map, one can easily check that

a) $f$ maps homeomorphically each face to one of the hemispheres
in  $\cP^1\setminus \rP^1$;

b) $f$ maps homeomorphically the boundary of each face of the net to
   $\rP^1$ (see e.g.,~\cite {Vi}).

Conditions a) and b) imply that
   $X(f)$ is properly labeled and oriented.
   \enddemo

   The following statement accomplishes the proof of Theorem 1.

\proclaim{Realization Theorem}
Let $X$ be a properly oriented and labeled planar chord diagram of
order $2n-2$,
   $C$  be a finite subset on $\rP^1$ consisting of $2n-2$ distinct points
   with a cyclic order induced
   by the positive orientation of  $\rP^1$. Label points in $C$ by the
   numbers
   $1\prec 2\prec \ldots \prec 2n-2\prec 1$ preserving the
    cyclic order.
   Then there exists a  degree $n$ generic $M$-function $f$ with $C(f)=C$
and $X(f)=X$,
whose properly oriented and labeled chord diagram coincides with
$X$.
\endproclaim

The proof of the Realization Theorem is completely parallel to the
proofs of similar results in~\cite {Vi, SV}. We provide a sketch below.

In order to verify that for any properly oriented and labeled chord diagram $X$
satisfying  1) and 2) there exists a rational function with the equivalent
net, we construct, using $X$, a branched covering satisfying
conditions a) and b) in the proof above.
Having done that, we induce the complex structure on the preimage from
that on the image, and the topological branched covering
transforms into a holomorphic one.
The construction of the topological covering goes as follows. Note
that it suffices to define it on the disc containing the chord
diagram, and then to glue the map $\cP^1$ to $\cP^1$ by pasting together
the initial map and its complex conjugate along the real cycle.
In order to determine our map on the disc, we first map the vertices of
the chord diagram to the corresponding critical values.
Then we map the chords connecting the vertices to the arcs of $\rP^1$
between the corresponding critical values, taking into account their
orientation.
Finally, the faces of the chord diagram are mapped homeomorphically
to the upper or lower hemispheres according to their orientations.
\qed
\medskip

\demo{Proof of Theorem 1'} The proof follows immediately from the proposition
below. 
Let $\pi(D)$ denote the number of all distinct plantings of the trees $T(D)$
and $T'(D)$ obtained from $T(D)$ by a reflection.

\proclaim{Proposition} The number of all distinct proper labelings of 
$D^+$ and $D^-$ equals $\pi(D)\sharp_r(D)$.
\endproclaim

\demo{Proof}
All the possible proper labelings of $D^+$ are obtained from 
the $\sharp_r(D)$ labelings with the root labeled by 1 via
the action of the cyclic group 
$C_{2e}$ on the labels, where $e$ is the number of chords in $D$. 
Moreover, each proper labeling is obtained in this
way exactly $|\Aut(D^+)|$ times. If $D$ is orientation--symmetric, 
then the contribution
of $D^-$ in the total number of nonequivalent labelings equals 0, since
any labeling of $D^-$ is equivalent to some labeling of $D^+$. If $D$ is
not orientation--symmetric, then the contribution
of $D^-$ equals the contribution of $D^+$.

On the other hand, there are exactly $2e$ possibilities to plant
the tree $T(D)$, and each planting is obtained exactly $|\Aut(T(D))|$ times.
Hence, $\pi(D)$ equals $2e/|\Aut(T(D))|$ if there exists an automorphism
of $T(D)$ and $T'(D)$ (preserving the orientation of $\bR^2$), and twice
this number otherwise.

Introduce
$$\align
\alpha(D)&=\left\{\alignedat2 &2|\Aut(D^+)|\quad &&
\text{if $D$ is orientation--symmetric,}\\
&|\Aut(D^+)| \quad &&\text{otherwise},
\endalignedat\right.\\
\beta(D)&=\left\{\alignedat2 &2|\Aut(T(D))|\quad &&
\text{if there exists an automorphism of $T(D)$ and $T'(D)$,}\\
&|\Aut(T(D))| \quad &&\text{otherwise}.
\endalignedat\right.
\endalign
$$
It remains to prove that $\alpha(D)=\beta(D)$.
The first of the above numbers can be interpreted as the number of 
automorphisms $\varphi$ of the pair $\{D^+,D^-\}$ such that $\varphi(D^+)$
and $\varphi(D^-)$ are two opposite orientations. The second of the
above numbers can be interpreted as the number 
of automorphisms $\psi$ of the pair $\{T(D),T'(D)\}$ such that $\psi(T(D))$
and $\psi(T'(D))$ differ by reflection. It is easy to see that these two
numbers coincide.
\qed\enddemo
\enddemo

\demo{\it Proof of Theorem 2} The crucial ingredient of this proof is the
   following modification of the well-known construction sending a
   branched covering to the set of its branching points, cf. e.g.
   \cite {A1} or \cite {L}.
Denote by  $\De_m$ the connected set of all unordered $m$-tuples of
pairwise distinct points on $\rP^1$. Obviously, the space $\De_{m}$
is fibered over $\rP^1$ with $(m-1)$-dimensional contractible fibers,
and therefore $\pi_{1}(\De_{m})=\Bbb Z$. (In order to convince
yourself think of $m$ electrons on $\rP^1$.) Note that the group
$PSL_{2}(\bR)$ of real linear fractional transformations $x\mapsto\frac
{ax+b}{cx+d}$ with $ad-bc=1$ acts freely on the space $\MRat_{n}$ of
generic rational M-functions and preserves all the critical values.
We now define the
    {\it generalized Lyashko-Looijenga mapping}
$$
    LL\: \MRat_n/PSL_{2}(\bR) \to \De_{2n-2}
$$
sending the equivalence class of a generic M-function to the unordered
set of its branching points.
Note that $LL$ is a map between the spaces of the same dimension,
and by Theorem 1 the number of inverse
images of any point in  $\De_{2n-2}$ is  constant.

\proclaim{Lemma 1}
The map  $LL\: \MRat_n/PSL_{2}(\bR)
\to \De_{2n-2}$ is a finite covering of degree
    $\sharp^{\bR}_{n}$. Moreover the action of $\pi_{1}(\De_{2n-2})$
on the fiber
    coincides with the action of the group
    $C_{2n-2}$ of the cyclic shifts of labels on the set of all properly
    directed and labeled plane chord diagrams of order $2n-2$.
\endproclaim

\demo{Proof}
To get a covering we need to show that

i) the map ${LL}$ is open, i.e. it is a local homeomorphism
onto the image;

ii) the map  $ {LL}$ is surjective.

Property i) is shown
in \cite {K} in a more general situation.
The above Realization Theorem provides the surjectivity.
Finally, Theorem 1 shows that the number of points in each fiber
equals $\sharp^{\bR}_{n}$. Let us study the action of $\pi_{1}(\De_{2n-2})$
on the fiber of $LL$. Without loss of generality we can, for example,
choose a base point $pt$ in $\De_{2n-2}$ presenting the
configuration of critical values forming some regular $(2n-2)$-gon $\Pi$ on
$\rP^1$.   The generator of  $\pi_{1}(\De_{2n-2})$ can be then
geometrically realized as the family $\Pi(t),\,t\in
[0,\frac {2\pi}{2n-2})$ of all regular $(2n-2)$-gons. Let us choose
some rational function $f$ with the set of critical values given by
$\Pi=\Pi(0)$.
Let us now rotate $\Pi$ in the family $\Pi(t)$ and follow the branch of
$f(t),\,f(0)=f$. All topological characteristics of the net
$S(f(t))$ (and therefore of $X(f(t))$, see proof of Theorem 1) will
be preserved. When $t$ reaches $\frac {2\pi}{2n-2}$ (i.e when we come
back to the same base point $pt$) the labels experience a cyclic shift by
$1$.
\qed
\enddemo

Theorem 2 is proven.
\enddemo

\heading \S 3. Calculating the number of proper labelings for a given
planar chord diagram \endheading

\demo{Proof of Theorem 3} Consider a properly oriented
planar chord diagram of $P_n$-type. Such a diagram contains two pairs
of {\it corners}, that is, two pairs of vertices connected both by a
chord and by a boundary arc (e.g. pairs $1,2$ and $3,4$ on
Fig.~1). Let us fix a corner $c$ in such a way that the chord
joining $c$ with the adjacent corner in the corresponding pair,
$c'$, is directed from $c$ to $c'$. We consider all the proper labelings of
the diagram satisfying the additional condition: the corner $c$ is
labeled by $1$.

Let us establish a bijection $\pi$ between the set of all such labelings and
the set of all 2up-2down permutations of length $2n-1$. First, we
start at $c$ and move along the chord to $c'$. The only
way to leave $c'$ is along the boundary arc (the one not leading back
to $c$). We proceed in this way, choosing each time a chord or a
boundary arc leading to a vertex not visited previously. Observe that
each time there exists a unique possibility to proceed, until we
visit two other corners, and the process stops. We thus get a unique
linear order on the vertex set of the diagram.
Let $\L=(l_1=1,l_2,\dots,l_{2n})$ be the sequence of labels of the
vertices written in this order.

Recall that all the labels together form a cyclically ordered set
$\Cal S=\{1\prec 2\prec\dots\prec 2n\prec 1\}$. To find the
permutation $\sigma=\pi(\L)$, it is convenient to consider its entries
as elements of a linearly ordered set $\Sigma=\{1<2<\dots<2n-1\}$. The
sequence $\L$ can be decomposed into (intersecting) cyclically ordered
quadruples $\{l_1\prec l_2\prec l_3\prec l_4\prec l_1\}$,
$\{l_3\prec l_4\prec l_5\prec l_6\prec l_3\}, \dots$, each
representing a face of the chord diagram. On the other hand, any
2up-2down permutation $\sigma=(\sigma_1,\dots,\sigma_{2n-1})$ can be
decomposed into (intersecting) linearly ordered triples
$\{\sigma_1<\sigma_2<\sigma_3\}$,
$\{\sigma_3>\sigma_4>\sigma_5\},\dots$. Bijection $\pi$ takes the
$k$th quadruple in $\L$ to the $k$th triple in $\sigma$. The main
building blocks of the bijection $\pi$ are the following two
bijections between quadruples and triples.

Let $\Cal R=\{r_1\prec r_2\prec\dots\prec r_m\prec r_1\}$ be a
cyclically ordered set, and let $r$ be an arbitrary element in $\Cal
R$. Consider the set $\Cal R^4_r$ of all quadruples of the form
$\{r\prec a\prec b\prec c\prec r\}$, $a,b,c\in\Cal R$, and denote by
$f_1$, $f_2$, $f_3$, $f_4$ the numbers of elements in $\Cal R$ lying
between $r$ and $a$, $a$ and $b$, $b$ and $c$, and $c$ and $r$,
respectively (see Figure~4). On the other hand, let $\Cal
T=\{t_1<t_2<\dots <t_{m-1}\}$ be a linearly ordered set, and let $\Cal
T_+^3$ and $\Cal T_-^3$ be the sets of all increasing and decreasing
triples in $\Cal T$, respectively. Transformation $\pi^+_{r,m}$ takes
the quadruple $\{r\prec a\prec b\prec c\prec r\}$ to the triple
$\{\tau_1<\tau_2<\tau_3\}$ such that the number of elements in $\Cal T$
less than $\tau_1$ equals $f_1$, the number of elements in $\Cal T$
lying between $\tau_1$ and $\tau_2$ equals $f_2$, and the number of
elements in $\Cal T$ lying between $\tau_2$ and $\tau_3$ equals $f_4$
(see Fig.~4).

   \vskip 15pt
\centerline{\hbox{\epsfysize=4cm\epsfbox{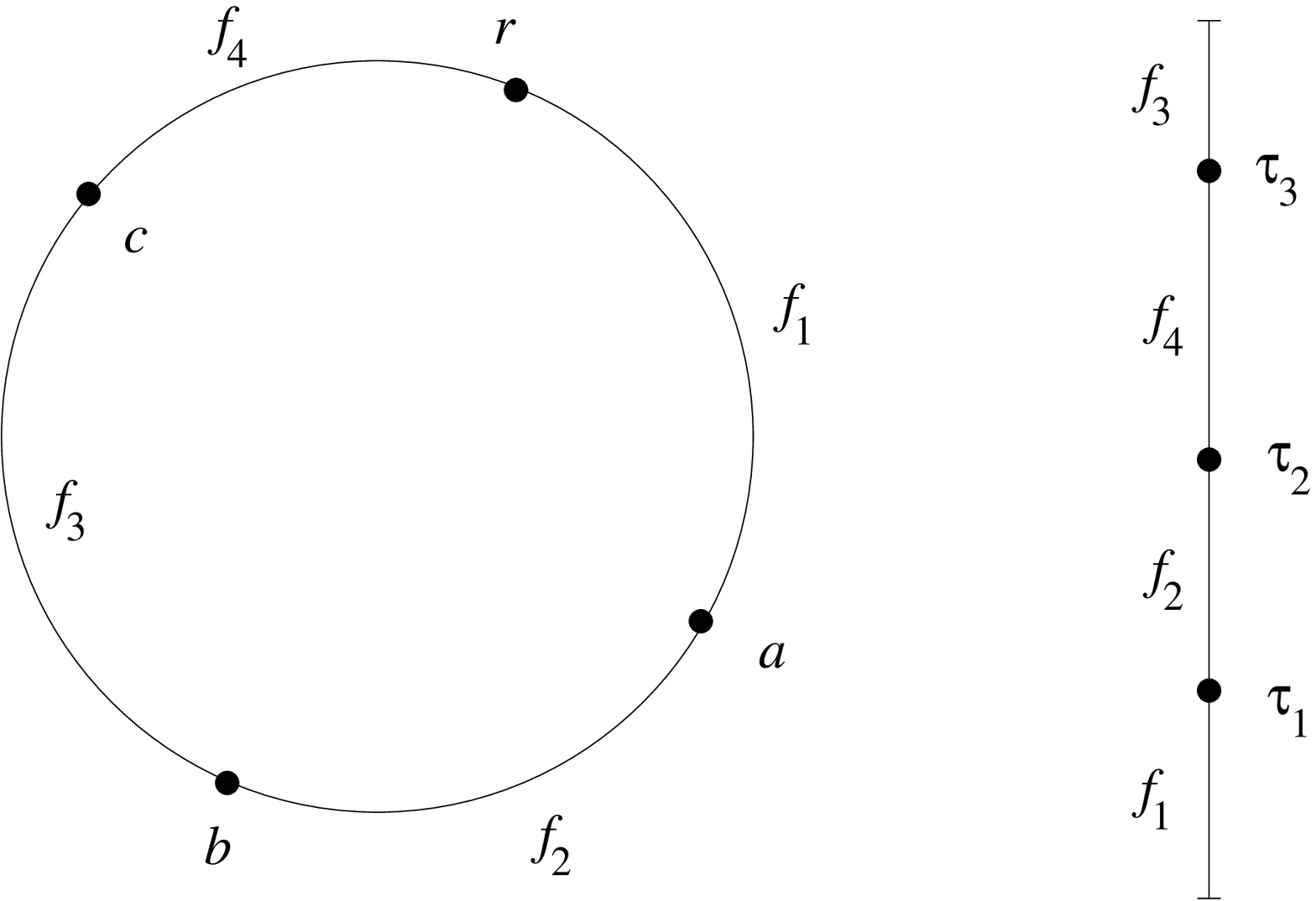}}}
\midspace{0.1mm} \caption{Fig.~4. Bijection $\pi^+_{r,m}$}
\medskip

It is easy to see that $\pi^+_{r,m}$ is a bijection between $\Cal
R^4_r$ and $\Cal T^3_+$, and that the number of elements in $\Cal T$
greater than $\tau_3$ equals $f_3$. Similarly,  $\pi^-_{r,m}$ takes
$\{r\prec a\prec b\prec c\prec r\}$ to the triple
$\{\tau'_1>\tau'_2>\tau'_3\}$ such that the number of elements in $\Cal T$
greater than $\tau'_1$ equals $f_1$, the number of elements in $\Cal T$
lying between $\tau'_1$ and $\tau'_2$ equals $f_2$, and the number of
elements in $\Cal T$ lying between $\tau'_2$ and $\tau'_3$ equals $f_4$.
It is easy to see that $\pi^-_{r,m}$ is a bijection between $\Cal
R^4_r$ and $\Cal T^3_-$, and that the number of elements in $\Cal T$
less than $\tau'_3$ equals $f_3$.

To establish the bijection $\pi$, we start from the quadruple
$\{1\prec l_2\prec l_3\prec l_4\prec 1\}$ and use $\pi^+_{1,2n}$ to
get $\{\sigma_1,\sigma_2,\sigma_3\}$. We now delete 1 and $l_2$ from
$\Cal S$, delete $\sigma_1$ and $\sigma_2$ from $\Sigma$, and use
$\pi^-_{l_3,2n-2}$ to get $\{\sigma'_3,\sigma'_4,\sigma'_5\}$. Observe
that the properties of $\pi^+_{r,m}$ and  $\pi^-_{r,m}$ stated above
ensure that $\sigma'_3=\sigma_3$, and hence we can glue the obtained
triples into a 2up-2down sequence of length 5. On the next step we
delete $l_3$ and $l_4$ from $\Cal S$, as well as $\sigma_3$ and
$\sigma_4$ from $\Sigma$, and use $\pi^+_{l_5,2n-4}$, and so on. On
the $k$th step we define $\Cal R=\Cal
S\setminus\{l_1,l_2,\dots,l_{2k-2}\}$, $\Cal
T=\Sigma\setminus\{\sigma_1,\sigma_2,\dots ,\sigma_{2k-2}\}$, $m=2n-2k+2$,
$\varepsilon=+$ if $k$ is even and $\varepsilon=-$ otherwise, and use
$\pi^\varepsilon_{l_{2k-1},m}$ to transform $\{l_{2k-1}\prec
l_{2k}\prec l_{2k+1}\prec l_{2k+2}\prec l_{2k-1}\}\in\Cal
R^4_{l_{2k-1}}$ into
$\{\sigma_{2k-1}<\sigma_{2k}<\sigma_{2k+1}\}\in\Cal T_\varepsilon^3$.

The rest of the proof follows immediately from the Proposition.
\qed
\enddemo

\demo{Proof of Corollary~1} We use the following result
  proved in \cite{CaS2, Theorem~1}. Let $A(k)$ denote the number of 2up-2down
permutations of length $k$; as mentioned in \S1, such permutations are
defined only for odd $k$. Then
$$\align
\sum_{n=0}^\infty
A(4n+1)\frac{x^{4n+1}}{(4n+1)!}&=\frac{\fy_0(x)\fy_1(x)-\fy_2(x)\fy_3(x)}
{\fy_0^2(x)-\fy_1(x)\fy_3(x)},\\
\sum_{n=0}^\infty
A(4n+3)\frac{x^{4n+3}}{(4n+3)!}&=\frac{\fy_3(x)}
{\fy_0^2(x)-\fy_1(x)\fy_3(x)},
\endalign
$$
where $\fy_j(x)$ is the $j$th Olivier function (see \S1 for the
definition). Corollary~1 follows immediately from this result and
Theorem~3. \qed\enddemo

\demo {Proof of Theorem 4} The proof goes along the same lines as the
proof of Theorem~3. The main difference is that the sequence of labels
is decomposed into intersecting cyclically ordered subsequences of sizes
$2d_1$, $2d_2,\dots$. As before, the intersection of two consequent
subsequences consists of two adjacent labels. In the
$P_k$-case these two labels formed a suffix of the former subsequence
and a prefix of the following one. In the present case, we still
maintain that the intersection is a prefix of the following subsequence;
however, it may lie somewhere in the middle of the former subsequence.

Accordingly, we decompose the permutation
$\sigma$ into linearly ordered subsequences of sizes $2d_1-1$,
$2d_2-1,\dots$, such that each odd subsequence is increasing and each even
subsequence is decreasing. As before, we establish bijections
between the  corresponding cyclical an linear subsequences. Assume
that we have a cyclically ordered $k$-tuple $\{a_1\prec
a_2\prec\dots\prec a_k\prec a_1\}$, and that the following cyclically
ordered subsequence intersects with this $k$-tuple by the elements
$a_j$ and $a_{j+1}$; evidently, $3\ls j\ls k-1$. Let $f_i$ denote the
number of elements in $\Cal
R$ lying between   $a_i$ and $a_{i+1}$, $i=1,\dots, k$. The we map
this $k$-tuple to the sequence $\{\tau_1<\tau_2<\dots>\tau_{k-1}\}$ such
that the number of elements in $\Cal T$ less than $\tau_1$ equals
$f_1$, the number of elements in $\Cal T$ greater than $\tau_{k-1}$
equals $f_j$, and the number of elements in $\Cal T$ lying between
$\tau_i$ and $\tau_{i+1}$ equals $f_{i+1}$ for $1\ls i\ls j-2$ and
$f_{i+2}$ for $j-1\ls i\ls k-2$.

The rest of the proof carries over without substantial changes.\qed
\enddemo

\end